\documentclass[11pt]{amsart}
\usepackage{amsmath, amssymb, mathptmx, mathtools, a4wide}
\usepackage{paralist}
\usepackage[utf8]{inputenc}
\usepackage[T1]{fontenc}

\numberwithin{equation}{section}

\theoremstyle{plain}
\newtheorem{thm}{Theorem}

\newtheorem*{prop*}{Proposition}
\newtheorem{lem}{Lemma}[section]

\newtheorem{cor}[lem]{Corollary}
\newcommand{\thmref}[1]{Theorem~\ref{#1}}
\newcommand{\lemref}[1]{lemma~\ref{#1}}

\theoremstyle{definition}

\newcommand{\mbf}{\mathbf}
\newcommand{\x}{\textbf}
\newcommand{\mf}{\mathbf}

\newcommand{\qq}{\qquad}

\newcommand{\mc}{\mathcal}

\newcommand{\mrm}{\mathrm}

\newcommand{\spz}{\mrm{Sp}_2(\mbf Z)}
\newcommand{\slz}{\mrm{SL}_2(\mbf Z)}

\newenvironment{psm}
  {\left(\begin{smallmatrix}}
  {\end{smallmatrix}\right)}
\begin{document}
 
\title{Bounds for the Petersson norms of the pullbacks of Saito-Kurokawa lifts}

\author{Pramath Anamby}
\address{Department of Mathematics\\
Indian Institute of Science\\
Bangalore -- 560012, India.}
\email{pramatha@iisc.ac.in}

\author{Soumya Das}
\address{Department of Mathematics\\
Indian Institute of Science\\
Bangalore -- 560012, India.}
\email{soumya@iisc.ac.in}

\begin{abstract}
Using the amplification technique, we prove that `mass' of the pullback of the Saito-Kurokawa lift of a Hecke eigen form $g\in S_{2k}$ is bounded by $k^{1-\frac{1}{210}+\epsilon}$. This improves the previously known bound $k$ for this quantity.
\end{abstract}


\maketitle

\section{Introduction}\label{intro}
Let $\ell$ be an integer and $M_\ell^2$ denote the space of Siegel modular forms of weight $\ell$ and degree $2$ on $\spz (\subseteq M_4(\mbf Z))$ and by $S_\ell^2$ the subspace of cusp forms. These are holomorphic functions defined on the Siegel upperhalf space $\mf H_2$ which consists of complex symmetric matrices $Z\in M_4(\mf C)$ whose imaginary part is positive-definite. If we write such a $Z=\begin{psm}\tau&z\\z&\tau'\end{psm}$, then $F|_{z=0}:=F\big( \begin{psm}\tau&0\\0&\tau'\end{psm}\big)$ is a modular form in $\tau$ and $\tau'$ (see \cite{fr} for more details) with weight $\ell$, which we call the pullback of $F$ to $\mf H \times \mf H$.

The study of pullbacks of automorphic forms has a rich history, see eg., \cite{blomer-khan-young}, \cite{ic}, \cite{ike}, \cite{liu-young}. In the context of Siegel modular forms, there are conjectures of Ikeda \cite{ike} relating such pullabacks to central values of $L$-functions. As an example, the Gross-Prasad conjecture would relate pullbacks of Siegel cusp forms of degree $2$ to central critical values of $L$-functions for $\mrm{GSp}(2) \times \mrm{GL}(2) \times \mrm{GL}(2)$. Ichino's beautiful result \cite{ic} studies this question for the Saito-Kurokawa (SK from now on) lifts of elliptic modular forms. Following the above notation, let us write $F|_{z=0} = \sum_{g_1,g_2} c_{g_1,g_2} g_1(\tau) g_2(\tau')$ (see also \cite{liu-young}), where $g_j$ runs over a Hecke basis of $S_\ell$, the space of elliptic cusp forms on $\slz$. Then Ichino proves that if $F=F_g$ is the SK-lift of $g\in S_{2 \ell -2}$ in the above, only the diagonal survives and $|c_{g_1,g_1}|$ is essentially given by the central value $L(1/2, \mrm{sym}^2g_1 \times g)$. 

It was moreover observed in \cite{liu-young} that comparison of the (normalised) norm of $F_g$ with the norm of its pullback provides a measure of the non-vanishing of the latter on average over the `projection' of $F_g|_{z=0}$ along $g_1 \times g_1$, as $g_1 \in S_{2\ell-2}$ varies. By a formula (see \eqref{nfgL}) in \cite{liu-young}, this also provides a measure of density of $F$ along $F|_{z=0}$ (see \cite[(1.13)]{liu-young}). This is made more precise in the next paragraph. We now make a change of notation, and use $2k$ for the weight $2 \ell -2$, in conformity to the afore-mentioned papers on the topic.

For an odd integer $k>0$, let $g\in S_{2k}$ be a normalized Hecke eigenform for $\slz$. Let $h\in S_{k+1/2}^{+}(\Gamma_0(4))$ be the Hecke eigenform associated to $g$ by the Shimura correspondence. Denote the Saito-Kurokawa lift of $g$ by $F_g\in S^2_{k+1}$. Let us define the quantity
\begin{equation}\label{def:nfg}
N(F_g):= \tfrac{1}{v_1^2}\left<F_g|_{z=0},F_g|_{z=0}\right> \big/  \tfrac{1}{v_2}\left<F_g,F_g\right>,
\end{equation}
where $v_1=\mrm{vol.}(\slz\backslash \mbf H)$ and $v_2=\mrm{vol.}(\spz\backslash \mbf H_2)$. Here $\left<F_g|_{z=0},F_g|_{z=0}\right>$ denotes the Petersson norm of $F_g|_{z=0}$ on $\slz\backslash \mbf H\times \slz\backslash \mbf H$ (see section \ref{sec:prelim} for more details).

Let $B_{k+1}$ denote the Hecke basis for $S_{k+1}$. Now Ichino's formula \cite{ic} immediately implies the following, as computed in \cite{liu-young}:
\begin{equation} \label{nfgL}
N(F_g)=\frac{\pi^2}{15} \left(  L(3/2,g)L(1,\mrm{sym}^2g) \right)^{-1} \cdot \frac{12}{k}\sum_{f\in B_{k+1}}
L(\frac{1}{2},\mrm{sym}^2f\times g).
\end{equation}
It is this quantity $N(F_g)$ that we are concerned with in this paper and we refer this quantity as the \textit{mass} of the pullback of the Saito-Kurokawa lift. Let us recall that as a special case of conjectures of \cite{CFKRS}, Liu and Young in \cite{liu-young} conjectured that $N(F_g)\sim 2$ as $k\rightarrow \infty$, and proved it on average over the family $g\in B_{2k}$ and $K\le k\le 2K$. In \cite{blomer-khan-young} a stronger asymptotic formula was obtained by considering only the smaller family $g\in B_{2k}$. Their result says that there exists some $\eta>0$ such that
\begin{equation}\label{eq:asymptotic}
\frac{12}{2k-1}\sum_{g\in B_{2k}}N(F_g)=2+O(k^{-\eta}).
\end{equation}
Now dropping all but one term, this asymptotic formula immediately gives $N(F_g)\ll k$. This bound is slightly better than the bound $N(F_g)\ll k\log k$  (cf. \cite{liu-young}) that one gets by using the convexity bound for $L(1/2, \mrm{sym}^2f\times g)$. In this paper we use the method of amplification (see \cite{iwnc-sarnak} for more details) to get a power-saving bound for an individual $N(F_g)$. We prove that
\begin{thm}\label{th:main}
For any $g\in B_{2k}$ and $\epsilon>0$, one has
\begin{equation}\label{eq:power}
N(F_g)\ll k^{1-\frac{1}{210}+\epsilon},
\end{equation}
where the implied constant depends only on $\epsilon$.
\end{thm}
Combining the asymptotic (\ref{eq:asymptotic}) and the power-saving obtained from \thmref{th:main} we have the following corollary, showing the existence of SK-lifts with non-vanishing `mass'.
\begin{cor}
We have
\begin{equation}
\# \{g\in B_{2k}\mid N(F_g)\neq 0\}\gg k^{1/210}.
\end{equation}
\end{cor}

As for the proof of \thmref{th:main}, we use the classical amplifier of Iwaniec-Sarnak \cite{iwnc-sarnak}. Now instead of inserting the amplifier in the sum  $\sum_{g\in B_{2k}}N(F_g)$, we insert the amplifier in a modified sum $\sum_{g\in B_{2k}}S_g$, where $S_g=L(3/2,g) N(F_g)$ (cf. \eqref{nfgL}). Since $L(3/2,g)^{-1}\ll 1$, this doesn't have any affect on the final bound for $N(F_g)$. Moreover, this modification helps in reducing the complexities of the further calculations. The proof of \thmref{th:main} follows a slightly different trajectory than that of the proof of the asymptotic \ref{eq:asymptotic}. In particular we have to be careful in keeping track of the dependence on the weight $2k$ throughout in a quantitative fashion.

Finally let us mention that one way of obtaining a better power saving seems to be an improvement in the error term that we obtain when we express $L(1,\mrm{sym}^2f)$ $(f\in B_k)$ as a Dirichlet polynomial (see lemma \ref{lem:dirpol}), which in turns relies on the results of Y.-K. Lau and J. Wu \cite{lau-wu}.

\noindent\textbf{Acknowledgements.} The first author is a DST- INSPIRE fellow at IISc, Bangalore and acknowledges the financial support from DST (India). The second author acknowledges financial support in parts from the UGC Centre for Advanced Studies, DST (India) and IISc, Bangalore during the completion of this work.

\section{Preliminaries}\label{sec:prelim}
In this section we collect some necessary results and formulae that will be used later in the paper. Throughout the article we follow the convention that $f(x)\ll g(x)$ ($g(x)\ge 0$) means there exist constants $M$ and $N$ such that $|f(x)|\le M\cdot g(x)$ for $x>N$. Moreover, $\epsilon$ will always denote an arbitrarily small positive
constant, but not necessarily the same one from one occurrence to
the next.
\subsection{The norm $\left<F_g|_{z=0},F_g|_{z=0}\right>$:} Let $F_g\in S^2_{k+1}$ be defined as in section \ref{intro}. Then $\left<F_g|_{z=0},F_g|_{z=0}\right>$ as in the definition of $N(F_g)$ in (\ref{def:nfg}) is the product of Petersson inner products on $\slz\backslash \mbf H\times \slz\backslash \mbf H$ and it is given by
\begin{equation}
\left<F_g|_{z=0},F_g|_{z=0}\right>=\int_{\slz\backslash \mbf H}\int_{\slz\backslash \mbf H}\mid F_g\big( \begin{psm}\tau& \\ &\tau'\end{psm}\big)\mid^2 \mrm{Im}(\tau)^{k+1}\mrm{Im}(\tau')^{k+1}d\mu(\tau)d\mu(\tau'),
\end{equation}
where $d\mu(z)=y^{-2}dxdy$, if $z=x+iy, y>0$.

\subsection{Petersson trace formula:} Let $f\in B_k$ and $\lambda_f(n)$ be the normalized Fourier coefficients of $f$. Then we have
\begin{equation}\label{eq:ptf}
\frac{2\pi^2}{k-1}\sum_{f\in B_k}\frac{\lambda_f(m)\lambda_f(n)}{L(1,\mrm{sym}^2f)}=\delta_{mn}+2\pi i^{-k}\sum_{c=1}^{\infty}\frac{S(m,n;c)}{c}J_{k-1}\left(\frac{4\pi\sqrt{mn}}{c}\right),
\end{equation}
where $\delta_{mn}=1$ if $m=n$ and $0$ otherwise, $S(m,n;c)$ is the Kloosterman sum and $J_{k-1}$ is the Bessel function. For the Bessel function we have the best possible upper bounds given by (see \cite{landau})
\begin{equation}\label{eq:bessel}
|J_k(x)|\ll \min\{k^{-1/3},|x|^{-1/3}\}
\end{equation}
for any real $x$ and $k\ge 0$. We use the following bound for the rapid decay of the Bessel function near zero.
\begin{equation}\label{eq:besseldecay}
J_k(x)\ll  \frac{|x/2|^k}{\Gamma(k+1)},\qq\text{ for } x>0.
\end{equation}
Let $D_1$ denote the off-diagonal term in (\ref{eq:ptf}). Then using (\ref{eq:besseldecay}) as above, one can truncate the $c$ sum in $D_1$ at $c\le 100 \frac{4\pi\sqrt{mn}}{k}$ upto a very small error, say $k^{-100}$. Then using (\ref{eq:bessel}) we can write
\begin{equation}\label{eq:bndoffd}
D_1=2\pi i^{-k}\sum_{c\le 100 \frac{4\pi\sqrt{mn}}{k}}\frac{S(m,n;c)}{c}J_{k-1}\left(\frac{4\pi\sqrt{mn}}{c}\right)+O(k^{-100})\ll \frac{\sqrt{mn}}{k^{4/3}}.
\end{equation}
Let $f\in S_k$ be a Hecke eigen eigenform of weight $k$ and $A_f(m,n)$ the Fourier-Whittaker coefficients of the symmetric square lift of $f$ (see \cite{gel-jac}). Then $A_f(m,n)$ is given by
\begin{equation}\label{eq:afmn}
A_f(m,n)=\sum_{d|(m,n)}\mu(d)A_f(m/d,1)A_f(n/d,1), \qquad \text{where} \qquad
A_f(r,1)=\sum_{ab^2=r}\lambda_f(a^2).
\end{equation}

\subsection{The approximate functional equation:} For $f\in B_{k+1}$ and $g\in B_{2k}$ we would use the approximate functional equation (cf. \cite{blomer-khan-young, iwnc-kw})
\begin{equation}\label{eq:aprxfe}
L(1/2,\mrm{sym}^2f\times g)=2\sum_{n,m\ge 1}\frac{\lambda_g(n)A_f(m,n)}{n^{1/2}m}W(nm^2).
\end{equation}
Here $W$ is a rapidly decaying smooth weight function. In fact $W$ satisfies (see \cite{blomer-khan-young})
\begin{equation}\label{eq:decayW}
x^{(j)}W^{(j)}(x)\ll_{j,A}\big(1+\frac{x}{k^2}\big)^{-A}
\end{equation}
for any $j,A\ge 0$. For our case we take
\begin{equation}\label{eq:W}
W(x)=\frac{1}{2\pi i}\int_{(1)}\frac{\Lambda_k(\frac{1}{2}+s)}{\Lambda_k(\frac{1}{2})}\left(\cos\frac{\pi s}{10 A}\right)^{-60 A}x^{-s}\frac{ds}{s},
\end{equation}
where (by invoking \cite[p.~2626]{blomer-khan-young}, by replacing $k$ with $k+1$, $\kappa$ with $k$, so that we are in the case $\kappa<k$)
\begin{equation}
\Lambda_k(s)=(2\pi)^{-3s}\Gamma(s+2k-\frac{1}{2})\Gamma(s+k-\frac{1}{2})\Gamma(s+\frac{1}{2}).
\end{equation}
We also require the following set-up in our proof. As in \cite{blomer-khan-young} we define
\begin{equation}\label{eq:mf}
\mc M_f(r):=\frac{12}{2k-1}\sum_{g\in B_{2k}}\lambda_g(r)\frac{L(1/2,\mrm{sym}^2f\times g)}{L(1,\mrm{sym}^2g)}.
\end{equation}
Using Deligne's bound, positivity and Theorem 1.4 in \cite{blomer-khan-young} we have
\begin{equation}\label{eq:mfrbound}
\mc M_f(r)\ll r^\epsilon\mc M_f(1)\ll r^\epsilon k^\epsilon.
\end{equation}
Also applying Petersson formula (\ref{eq:ptf}) to (\ref{eq:mf}) we get $\mc M_f(r)=\mc M^{(1)}_f(r)+\mc M^{(2)}_f(r)$, where $\mc M^{(1)}_f(r)$ is the diagonal contribution and $\mc M^{(2)}_f(r)$ is the off-diagonal contribution. We have the following
\begin{equation}\label{eq:mf1}
\mc{M}^{(1)}_f(r)=\frac{2}{\zeta(2)}\sum_{m}\frac{A(m,r)}{r^{1/2}m}W(m^2).
\end{equation}
For the off-diagonal term we have
\begin{equation}
\mc M^{(2)}_f(r)=\frac{4\pi i^k}{\zeta(2)}\sum_{n,m,c\ge 1}\frac{A(n,m)}{m^{1/2}n}W(mn^2)\frac{S(m,r,c)}{c}J_{2k-1}\left(\frac{4\pi \sqrt{mr}}{c}\right).
\end{equation}
Now using the rapid decay of $W$ as in (\ref{eq:decayW}), we can truncate the $m-$sum at $m\le k^{2+\epsilon}n^{-2}$ upto a negligible error, say $k^{-100}$. If one wants to bound $\mc M^{(2)}_f(r)$, then by the insertion of smooth partitions of unity for the $m$ and $c$ sums, it is enough to bound the quantity
\begin{equation}
\mc{M}^{(2)}_f(r, M, C)=\sum_{n,c}\frac{\Omega_1(c/C)}{nCM^{1/2}}\sideset{}{^*}\sum_{d( c)} e(\frac{dr}{c})\sum_{m}A(n,m)e(\frac{\overline{d}m}{c})\Omega_2(\frac{m}{M})J_{2k-1}\big(\frac{4\pi\sqrt{mr}}{c}\big)
\end{equation}
for
\begin{equation}\label{eq:boundCM}
M\le \tfrac{k^{2+\epsilon}}{n^2},\qq C\le 100\tfrac{\sqrt{Mr}}{k}.
\end{equation}
The truncation over $c$ comes from the rapid decay of the Bessel function near $0$ (see (\ref{eq:besseldecay})). Here $\Omega_1$ and $\Omega_2$ are fixed, smooth, compactly supported weight functions. From (\ref{eq:boundCM}) we immediately get that
\begin{equation}\label{eq:boundMcm}
k^2r^{-1}\le M\le k^{2+\epsilon},\qq cn\ll r^{1/2}k^\epsilon.
\end{equation}
Also using the Voronoi formula we can write (see \cite{blomer-khan-young} and \cite{ms})
\begin{equation}
\mc{M}^{(2)}_f(r, M, C)=\sum_{n,c}\frac{\Omega_1(c/C)}{nCM^{1/2}}c\sum_{m_1|c}\sum_{\pm}\sum_{m_2}\frac{A(m_2,m_1)}{m_1m_2}\sideset{}{^*}\sum_{d( c)} e(\frac{dr}{c})S(nd,\pm m_2,c/m_1)\Psi^{\pm}\big(\frac{m_2m_1^2}{c^3n}\big).
\end{equation}
We also have (see lemma 5.1 in \cite{blomer-khan-young})
\begin{equation}\label{eq:boundPsi}
\Psi^{\pm}(x)\ll_{A,\varepsilon} k^\varepsilon \frac{xc^2}{r^{1/4}}\big(1+\frac{xc^3}{k^\varepsilon M^{1/2}r^{3/2}}\big)^{-A},\qq \text{ for }x\ge \frac{1}{c^3n}
\end{equation}
which follows from a very careful estimation of certain oscillatory integrals (see \cite{blomer-khan-young}) and
\begin{equation}\label{eq:boundksum}
|\sideset{}{^*}\sum_{d( c)} e(\frac{dr}{c})S(nd,\pm m_2,c/m_1)|\le c\tau (c)(c,n).
\end{equation}
We also make use of the useful observation that almost all $L(1,\mrm{sym}^2f)$ for $f\in B_{k+1}$ can be approximated by a convergent Dirichlet series with rapidly decaying weight function (see lemma 6.1 in \cite{blomer-khan-young}). More precisely:
\begin{lem}\label{lem:dirpol}
Given $\delta_1,\delta_2>0$, there is a $\delta_3$ such that
\begin{equation}\label{eq:dirpol}
L(1,\mrm{sym}^2f)=\sum_{d_1,d_2}\frac{\lambda_f(d_1^2)}{d_1d_2^2}\exp\left(-\frac{d_1d_2^2}{k^{\delta_1}}\right)+O(k^{-\delta_3})
\end{equation}
for all but $O(k^{\delta_2})$ cusp forms $f\in B_{k+1}$.
\end{lem}
We also note from the proof of \lemref{lem:dirpol} in \cite{blomer-khan-young} that one can take, $\delta_3<\delta_1\delta_2/62$.

\section{Proof of \thmref{th:main}}
Recall that $S_g=L(3/2,g)N(F_g)$. We start we the following sum with the amplifier given by $\mid\sum_{n\le N}\alpha_n\lambda_g(n)\mid^2$
\begin{equation}\label{eq:The sum}
\mbf S_A=\frac{12}{2k-1}\sum_{g\in B_{2k}}\mid\sum_{n\le N}\alpha_n\lambda_g(n)\mid^2S_g.
\end{equation}
Expanding the sum (\ref{eq:The sum}) and using the Hecke relation $\lambda_g(m)\lambda_g(n)=\underset{d|(m,n)}{\sum}\lambda_g(\frac{mn}{d^2})$, we have
\begin{equation}\label{eq:afterheckerelations}
\mbf S_A=\frac{12}{2k-1}\sum_{g\in B_{2k}}\sum_{n_1,n_2\le N}\alpha_{n_1}\overline{\alpha_{n_2}}\sum_{d|(n_1,n_2)}\lambda_g(\frac{n_1n_2}{d^2}) S_g.
\end{equation}
Now substituting for $S_g$ and using the approximate functional equation (see (\ref{eq:aprxfe})) and the Petersson formula (see (\ref{eq:ptf})) for the sum over $g$  we get
\begin{equation*}
\mbf S_A=\frac{\pi^2}{15}\cdot \frac{12}{k}\sum_{n_1,n_2\le N}\alpha_{n_1}\overline{\alpha_{n_2}}\sum_{f\in B_{k+1}}\sum_{d|(n_1,n_2)}\left(\mc{M}^{(1)}_f(\frac{n_1n_2}{d^2})+\mc{M}^{(2)}_f(\frac{n_1n_2}{d^2})\right),
\end{equation*}
where the quantities $\mc{M}^{(i)}(r)$ for $i=1,2$ are as in section \ref{sec:prelim}. To apply the Petersson formula for the sum over $f$, it is convenient to introduce the quantity $L(1,\mrm{sym}^2f)$ in the above sum. Thus we write
\begin{equation*}
\mbf S_A=\frac{\pi^2}{15}\cdot \frac{12}{k}\sum_{n_1,n_2\le N}\alpha_{n_1}\overline{\alpha_{n_2}}\sum_{f\in B_{k+1}}\frac{L(1,\mrm{sym}^2f)}{L(1,\mrm{sym}^2f)}\sum_{d|(n_1,n_2)}\left(\mc{M}^{(1)}_f(\frac{n_1n_2}{d^2})+\mc{M}^{(2)}_f(\frac{n_1n_2}{d^2})\right).
\end{equation*}
Making a change of variables and rearranging the summation we get
\begin{equation*}
\mbf S_A=\frac{\pi^2}{15}\cdot \frac{12}{k}\sum_{d\le N}\sum_{n_1,n_2\le N/d}\alpha_{dn_1}\overline{\alpha_{dn_2}}\sum_{f\in B_{k+1}}\frac{L(1,\mrm{sym}^2f)}{L(1,\mrm{sym}^2f)}\left(\mc{M}^{(1)}_f(n_1n_2)+\mc{M}^{(2)}_f(n_1n_2)\right).
\end{equation*}
Now we use \lemref{lem:dirpol} to get
\begin{equation}\label{eq:smain}
\begin{aligned}
\mbf S_A&=\frac{\pi^2}{15}\cdot \frac{12}{k}\sum_{d\le N}\sum_{n_1,n_2\le N/d}\alpha_{dn_1}\overline{\alpha_{dn_2}}\sum_{f\in B_{k+1}}\frac{1}{L(1,\mrm{sym}^2f)}\sum_{d_1,d_2}\frac{\lambda_f(d_1^2)}{d_1d_2^2}\exp\left(-\frac{d_1d_2^2}{k^{\delta_1}}\right)\\
&\times\left(\mc{M}^{(1)}_f(n_1n_2)+\mc{M}^{(2)}_f(n_1n_2)\right)+ O \left(\big(N^{\epsilon}\sum_{n_1,n_2\le N}|\alpha_{n_1}\overline{\alpha_{n_2}}|\big) \big(k^{-\delta_3+\epsilon}+k^{\delta_2-1+\epsilon}\big)\right).
\end{aligned}
\end{equation}
As in \cite{blomer-khan-young}, the error term comes from two sources: the error in \lemref{lem:dirpol} and the forms $f\in B_{k+1}$ for which (\ref{eq:dirpol}) doesn't hold. In the former case we use that the sum over $d_1,d_2$ is bounded by $\log(k)$ (from the bound of \cite{HL} and \lemref{lem:dirpol}); in the latter case we estimate trivially using (\ref{eq:mfrbound}).

Denote by $\mbf S_A^{(1)}$ and $\mbf S_A^{(2)}$ the terms corresponding to $\mc{M}^{(1)}$ and $\mc{M}^{(2)}$ in (\ref{eq:smain}) respectively. We first bound $\mbf S_A^{(1)}$.

\subsection{The term $\mbf S_A^{(1)}$}
Expanding out $\mc M_f^{(1)}$ from (\ref{eq:mf1}) we have
\begin{equation*}
\mc{M}^{(1)}_f(n_1n_2)=\frac{2}{\zeta(2)}\underset{a_2b_2^2|n_1n_2}{\sum_{a_1,a_2,b_1,b_2}}\sum_{d_4|(a_1^2,a_2^2)}\mu(\frac{n_1n_2}{a_2b_2^2})\frac{a_2b_2^2\lambda_f(\frac{a_1^2a_2^2}{d_4^2})}{(n_1n_2)^{3/2}a_1b_1^2}W(\frac{n_1^2n_2^2a_1^2b_1^4}{a_2^2b_2^4}).
\end{equation*}
\subsubsection{The diagonal:}
We apply the Petersson formula (\ref{eq:ptf}) for the sum over $f$ in $\mbf S_A^{(1)}$ and denote by $\mbf S_A^{(11)}$ the corresponding diagonal term. Then we have
\begin{equation*}
\mbf S_A^{(11)}=\frac{2\pi^2}{15\zeta(2)^2}\sum_{d\le N}\sum_{n_1,n_2\le N/d}\alpha_{dn_1}\overline{\alpha_{dn_2}}\underset{a_2b_2^2|n_1n_2}{\sum_{a_1,a_2,b_1,b_2,d_1,d_2}}\underset{d_1d_4=a_1a_2}{\sum_{d_4|(a_1^2,a_2^2)}}\tfrac{\mu(\frac{n_1n_2}{a_2b_2^2})a_2b_2^2}{(n_1n_2)^{3/2}a_1b_1^2d_1d_2^2}W(\tfrac{n_1^2n_2^2a_1^2b_1^4}{a_2^2b_2^4})\exp\left(-\tfrac{d_1d_2^2}{k^{\delta_1}}\right).
\end{equation*}
Using Mellin inversion we can write
\begin{equation*}
\mbf S_A^{(11)}=\frac{2\pi^2}{15\zeta(2)^2}\int_{(1)}\int_{(1)}\zeta(2+4u)\zeta(2+2v)\zeta(1+u+\frac{v}{2})\widetilde{W}(u)\Gamma(v)k^{\delta_1 v}B_N(u,v)\frac{du}{2\pi i}\frac{dv}{2\pi i},
\end{equation*}
where
\begin{equation*}
B_N(u,v)=\sum_{1\le d\le N}\sum_{n_1,n_2\le N/d}\alpha_{dn_1}\overline{\alpha_{dn_2}}\underset{a_2b_2^2|n_1n_2}{\sum_{a_2,b_2}}\frac{\mu(\frac{n_1n_2}{a_2b_2^2})\sigma_{v/2-u}(a^2)b_2^{2+4u}}{(n_1n_2)^{3/2+2u}a_2^{v-2u}}
\end{equation*}
and from \eqref{eq:W},
\begin{equation}
\widetilde{W}(u)=\frac{\Lambda_k(\frac{1}{2}+u)}{\Lambda_k(\frac{1}{2})}\left(\cos\frac{\pi u}{10 A}\right)^{-60 A}\cdot \frac{1}{u}.
\end{equation}
We further denote by $|B_N(u,v)|$ the sum as in $B_N(u,v)$, but with all terms replaced by their absolute values. First we move the line of integration w.r.t $u$ to $-\delta$, for some $1/2<\delta<1$ and encounter the pole of $\tilde{W}$ at $u=0$ and the poles of $\zeta$ at $u=-v/2$ and $u=-1/4$. Then the integral over $u$ equals
\begin{equation*}
\begin{aligned}
R(v):=&\zeta(2)\zeta(1+\frac{v}{2})B_N(0,v)+\zeta(2-2v)\widetilde{W}(-\frac{v}{2})B_N(-\frac{v}{2},v)+\zeta(\frac{3}{4}+\frac{v}{2})\widetilde{W}(-\frac{1}{4})B_N(-\frac{1}{4},v)\\
& \qquad+\int_{(-\delta)}\zeta(2+4u)\zeta(1+u+\frac{v}{2})\widetilde{W}(u)B_N(u,v)\frac{du}{2\pi i}.
\end{aligned}
\end{equation*}
Let us call by $R_1(v),R_2(v)$ the functions of $v$ appearing on the first and second line in the above expression for $R(v)$. Next we move the line of integration w.r.t $v$ to $\epsilon>0$ and cross the pole of $\zeta$ at $v=1/2$. The contribution to the residue only comes from the last two terms in $R_1(v)$. Then the contribution from $R_1$ to $\mbf S_A^{(11)}$ becomes (a sum of four terms):
\begin{equation}
R_1:= (|B_N(0,\epsilon)|+\mid B_N(-\frac{\epsilon}{2},\epsilon)\mid) k^{\epsilon}+\mid B_N(-\frac{1}{4},\frac{1}{2})\mid k^{\delta_1/2-1/2} +\mid B_N(-\frac{1}{4},\epsilon)\mid k^{-1/2+\epsilon} .
\end{equation} 
Since $R_2(v)$ is entire, its contribution to $\mbf S_A^{(11)}$ is just the integral over the two new lines of integrations, namely
\begin{equation}
R_2:= \mid B_N(-\delta, \epsilon) \mid k^{-2 \delta+\epsilon}.
\end{equation}
Noting the following bounds: $|B_N(0,\epsilon)|\le |B_N(0,0)|$, $\mid B_N(-\frac{1}{4},\frac{1}{2})\mid\le \mid B_N(-\frac{1}{4},\epsilon)\mid\le \mid B_N(-\frac{1}{4},0)\mid$ and $\mid B_N(-\delta, \epsilon) \mid\le \mid B_N(-\delta, 0) \mid$, we have
\begin{equation}\label{eq:sa11}
\mbf S_A^{(11)} \ll (|B_N(0,0)|+\mid B_N(-\frac{\epsilon}{2},\epsilon)\mid) k^{\epsilon}+\mid B_N(-\frac{1}{4},0)\mid (k^{\delta_1/2 -1/2})+\mid B_N(-\delta, 0) \mid k^{-2\delta+\epsilon}.
\end{equation}
\subsubsection{The off-diagonal:}
Denote the off-diagonal terms of $\mbf S_A^{(1)}$ by $\mbf S_A^{(12)}$, then we have
\begin{align*}
\mbf S_A^{(12)}&=2\pi i^{-k}\frac{2\pi^2}{15\zeta(2)^2}\sum_{d\le N}\sum_{n_1,n_2\le N/d}\alpha_{dn_1}\overline{\alpha_{dn_2}}\underset{a_2b_2^2|n_1n_2}{\sum_{a_1,a_2,b_1,b_2,d_1,d_2}}\sum_{d_4|(a_1^2,a_2^2)}\sum_{c}\frac{\mu(\frac{n_1n_2}{a_2b_2^2})a_2b_2^2}{c(n_1n_2)^{3/2}a_1b_1^2d_1d_2^2}\\
&\times S\big(\frac{a_1^2a_2^2}{d_4^2},d_1^2,c\big)W(\frac{n_1^2n_2^2a_1^2b_1^4}{a_2^2b_2^4})\exp\left(-\frac{d_1d_2^2}{k^{\delta_1}}\right)J_k\big(\frac{4\pi a_1a_2d_1}{cd_4}\big).
\end{align*}
We can truncate the sum over $c$ at $c\le 100\frac{4\pi a_1a_2d_1}{d_4 k}$ by using the rapid decay of Bessel function near $0$ (see (\ref{eq:besseldecay})). Next we use the bound (\ref{eq:decayW}) for $W$ with $j=0$ and $A=\frac{1}{2}+\frac{\epsilon}{2}$ and the trivial bounds $|S(*,*,c)|\le c$, $J_k(x)\ll k^{-1/3}$ (see (\ref{eq:bessel})) to  see that
\begin{equation}\label{eq:S12}
\mbf S_A^{(12)}\ll \left(\sum_{d\le N}\sum_{n_1,n_2\le N/d}|\alpha_{dn_1}\overline{\alpha_{dn_2}}|\underset{a_2b_2^2|n_1n_2}{\sum_{a_2,b_2}}\frac{\mid\mu(\frac{n_1n_2}{a_2b_2^2})\mid a_2^{2+\epsilon}b_2^{4+2\epsilon}}{(n_1n_2)^{5/2+\epsilon}}\right)k^{-1/3+\delta_1+\epsilon}.
\end{equation}

\subsection{The term $\mbf S_A^{(2)}$}
Now we proceed to estimate $\mbf S_A^{(2)}$. From the arguments in section \ref{sec:prelim} it is enough to bound
\begin{equation*}
\mbf S_A^{(2)}(M,C):=\frac{12}{k}\sum_{d\le N}\sum_{n_1,n_2\le N/d}\alpha_{dn_1}\overline{\alpha_{dn_2}}\sum_{f\in B_{k+1}}\frac{1}{L(1,\mrm{sym}^2f)}\sum_{d_1,d_2}\tfrac{\lambda_f(d_1^2)}{d_1d_2^2}\exp\left(-\tfrac{d_1d_2^2}{k^{\delta_1}}\right)\mc{M}^{(2)}_f(n_1n_2, M, C).
\end{equation*}
We use (\ref{eq:boundPsi}) and (\ref{eq:boundksum}) in the above equation and also note that the using the exponential decay, the sum over $d_1$ can be truncated at $d_1\le k^{\delta_1+\epsilon}$ with a very small error. Thus we are left to bound
\begin{equation}\label{eq:s2mc}
\mbf S_A^{(2)}(M,C)\ll \sum_{d\le N}\sum_{n_1,n_2\le N/d}|\alpha_{dn_1}\overline{\alpha_{dn_2}}|\left(\sum_{d_1\le k^{\delta_1+\epsilon}}\sum_{n}\sum_{C\le c\le 2C} T(d_1,n_1n_2,n,c,M)+O(k^{-100})\right),
\end{equation}
where
\begin{align*}
T(d_1,n_2n_1,n,c,M)&=k^\varepsilon\underset{m_1|c}{\sum_{m_2m_1^2\le k^\epsilon M^{1/2}(n_1n_2)^{3/2}n}}\tfrac{m_1\tau(c)(c,n)}{(n_1n_2)^{1/4}n^2d_1M^{1/2}}\left|\frac{12}{k}\sum_{f\in B_{k+1}}\tfrac{\lambda_f(d_1^2)A(m_2,m_1)}{L(1,\mrm{sym}^2f)}\right|\\
&\ll k^\varepsilon\underset{al_1^2m_1|c}{\underset{a^3l_1^4l_2^2m_2m_1^2\le k^\epsilon M^{1/2}(n_1n_2)^{3/2}n}{\sum_{a,l_1,l_2,m_1,m_2}}}\tfrac{al_1^2m_1\tau(c)(c,n)}{(n_1n_2)^{1/4}n^2d_1M^{1/2}}\sum_{h|(m_1^2,m_2^2)}\left|\frac{12}{k}\sum_{f\in B_{k+1}}\tfrac{\lambda_f(d_1^2)\lambda_f(m_1^2m_2^2/h^2)}{L(1,\mrm{sym}^2f)}\right|.
\end{align*}
We use (\ref{eq:afmn}) and the Hecke relations to arrive at the previous step.

Now we apply the Petersson formula and using the rapid decay of Bessel function near $0$ for the off-diagonal term (see (\ref{eq:bndoffd})) we get with the same conditions on the variables as above that
\begin{equation}
T(d_1,n_2n_1,n,c,M)\ll k^\varepsilon\sum_{a,l_1,l_2,m_1,m_2}\tfrac{al_1^2m_1\tau(c)(c,n)}{(n_1n_2)^{1/4}n^2d_1M^{1/2}}\sum_{h|(m_1^2,m_2^2)}\left(\delta_{d_1h=m_1m_2}+O\big(\frac{d_1m_1m_2}{hk^{4/3}}\big)\right).\nonumber
\end{equation}
Since the sum over $l_2$ is free, it is $\ll \left(\frac{k^\epsilon M^{1/2}(n_1n_2)^{3/2}n}{a^3l_1^4m_2m_1^2}\right)^{1/2}$. Thus
\begin{equation}
T(d_1,n_2n_1,n,c,M)\ll k^\varepsilon\underset{al_1^2m_1|c}{\underset{a^3l_1^4m_2m_1^2\le k^\epsilon M^{1/2}(n_1n_2)^{3/2}n}{\sum_{a,l_1,m_1,m_2}}}\tfrac{(n_1n_2)^{1/2}\tau(c)(c,n)}{a^{1/2}m_2^{1/2}n^{3/2}d_1M^{1/4}}\sum_{h|(m_1^2,m_2^2)}\left(\delta_{d_1h=m_1m_2}+O\big(\frac{d_1m_1m_2}{hk^{4/3}}\big)\right).\nonumber
\end{equation}
Now the sum over $l_1$ is $\ll \left(\frac{k^\epsilon M^{1/2}(n_1n_2)^{3/2}n}{a^3m_2m_1^2}\right)^{1/4}$ and the sum over $a$ is $\ll \left(\frac{k^\epsilon M^{1/2}(n_1n_2)^{3/2}n}{m_2m_1^2}\right)^{-1/12}$. Thus
\begin{equation*}
T(d_1,n_2n_1,n,c,M)\ll k^\varepsilon\tfrac{(n_1n_2)^{3/4}\tau(c)(c,m)}{n^{7/6}d_1M^{1/6}}\underset{m_1|c}{\underset{m_2m_1^2\le k^\epsilon M^{1/2}(n_1n_2)^{3/2}n}{\sum_{m_1,m_2}}}\tfrac{1}{m_2^{2/3}m_1^{1/3}}\sum_{h|(m_1^2,m_2^2)}\left(\delta_{d_1h=m_1m_2}+O\big(\tfrac{d_1m_1m_2}{hk^{4/3}}\big)\right).
\end{equation*}
Now we evaluate the two inside sums separately.
Let $T_1$ and $T_2$ correspond to the diagonal and off-diagonal terms respectively in the above sum. Then with the same conditions on the variables as above, we have
\begin{equation*}
T_1:=\sum_{m_1,m_2}\frac{1}{m_2^{2/3}m_1^{1/3}}\sum_{h}\delta_{d_1h=m_1m_2}.
\end{equation*}
Making the following change of variable, $m_2^2=d_1^2h^2/m_1^2$ we get
\begin{equation*}
T_1=\underset{m_1|c, h|m_1^2}{\underset{hm_1\le k^\epsilon M^{1/2}(n_1n_2)^{3/2}n/d_1}{\sum_{h,m_1}}}\frac{m_1^{1/3}}{h^{2/3}d_1^{2/3}}\ll \underset{m_1|c}{\underset{m_1\le k^\epsilon M^{1/2}(n_1n_2)^{3/2}n/d_1}{\sum_{m_1}}}\frac{m_1^{1/3+\epsilon}}{d_1^{2/3}}\ll \frac{\sigma_{1/3+\epsilon}(c)}{d_1^{2/3}}.
\end{equation*}
Now consider the second sum
\begin{equation*}
T_2:=\underset{m_1|c}{\sum_{m_1,m_2}}\frac{1}{m_2^{2/3}m_1^{1/3}}\sum_{h|(m_1^2,m_2^2)}\frac{d_1m_1m_2}{hk^{4/3}}=k^{-4/3}d_1\underset{m_1|c}{\sum_{m_1,m_2}}\sum_{h|(m_1^2,m_2^2)}\frac{m_1^{2/3}m_2^{1/3}}{h}.
\end{equation*}
In both of the above sums $m_2m_1^2\le k^\epsilon M^{1/2}(n_1n_2)^{3/2}n$. Put $m_1^2=hm_3$. Then
\begin{equation*}
T_2=k^{-4/3}d_1\underset{hm_3|c^2, h|m_2^2}{\underset{hm_3m_2\le k^\epsilon M^{1/2}(n_1n_2)^{3/2}n}{\sum_{h,m_3,m_2}}}\frac{m_3^{1/3}m_2^{1/3}}{h^{2/3}}.
\end{equation*}
Now using the fact that $hm_3|c^2$ we find that the sums over $h$ and $m_3$ are $\ll c^{\epsilon}\sigma_{1/3}(c^2)$. The remaining sum over $m_2$ is then $\ll k^\epsilon M^{2/3}(n_1n_2)^{2}n^{4/3}$. This implies $T_2\ll d_1k^{-4/3+\epsilon}c^{2\epsilon}\sigma_{1/3}(c^2)M^{2/3}(n_1n_2)^{2}n^{4/3}$. Thus we have
\begin{equation*}
T(d_1,n_2n_1,n,c,M)\ll k^\varepsilon\frac{(n_1n_2)^{3/4}\tau(c)(c,n)}{n^{7/6}d_1M^{1/6}}\left(\frac{\sigma_{1/3+\epsilon}(c)}{d_1^{2/3}}+d_1k^{-4/3+\epsilon}c^{2\epsilon}\sigma_{1/3}(c^2)M^{2/3}(n_1n_2)^{2}n^{4/3}\right).
\end{equation*}
Using that $\sigma_\alpha(c)\ll c^\alpha $, $(c,n)\le c$ and that $M\le k^{2+\epsilon}$ from (\ref{eq:boundMcm}), we get
\begin{equation}\label{eq:T}
T(d_1,n_1n_2,n,c,M)\ll (n_1n_2)^{11/4+\epsilon}c^{5/3+\epsilon}n^{1/6}k^{-1/3+\epsilon}.
\end{equation}
Finally, using (\ref{eq:T}) in (\ref{eq:s2mc}) we use $c\ll (n_1n_2)^{1/2}k^{\epsilon}n^{-1}$ (from (\ref{eq:boundMcm})) and note that sum over $n$ is $\ll (n_1n_2)^{-3/4}k^{\epsilon}$, we get
\begin{equation}\label{eq:s2}
\mbf S_A^{(2)}(M,C)\ll \left(N^{20/3+\epsilon}\sum_{n_1,n_2\le N}|\alpha_{n_1}\overline{\alpha_{n_2}}|\right)k^{-1/3+\delta_1+\epsilon}.
\end{equation}
Putting everything together from (\ref{eq:sa11}), (\ref{eq:S12}) and (\ref{eq:s2}), we have
\begin{equation}\label{eq:SAfinal}
\begin{aligned}
\mbf S_A&\ll (|B_N(0,0)|+\mid B_N(-\frac{\epsilon}{2},\epsilon)\mid) k^{\epsilon}+\mid B_N(-\frac{1}{4},0)\mid k^{\delta_1/2 -1/2}+\mid B_N(-\delta, 0) \mid k^{-2\delta+\epsilon}\\
&+\mbf S_A^{(12)}+\left(\sum_{n_1,n_2\le N}|\alpha_{n_1}\overline{\alpha_{n_2}}|\right) (N^{20/3+\epsilon}k^{-1/3+\delta_1+\epsilon}+N^{\epsilon}k^{-\delta_3+\epsilon}+N^{\epsilon}k^{\delta_2-1+\epsilon}).
\end{aligned}
\end{equation}
where $\delta_1,\delta_2>0$ are arbitrary, $1/2<\delta<1$ and $\delta_3<\delta_1\delta_2/62$.
\subsection{Choice of the amplifier: }For a fixed $g_0$ in the sum (\ref{eq:The sum}) we choose the $\alpha_n$s following Iwaniec-Sarnak (\cite{iwnc-sarnak}) as below
\begin{equation}\label{eq:alpha}
\alpha_n=\begin{cases}
\lambda_{g_0}(p), &\text{ if } n=p\le N^{1/2};\\
-1, &\text{ if }n=p^2\le N;\\
0,& \text{ otherwise }.
\end{cases}
\end{equation}
Substituting in (\ref{eq:The sum}) and using the Hecke relation $\lambda_{g_0}(p)^2-\lambda_{g_0}(p^2)=1$, we find that
\begin{equation}\label{eq:amplifiedterm}
\frac{12}{2k-1}S_{g_0}\big|\sum_{p\le N^{1/2}}1\big|^2\le \mbf S_A.
\end{equation}
We proceed to bound the quantities $|B_N(0,0)|$, $\mid B_N(-\delta, 0) \mid$, $|B_N(-\epsilon/2,\epsilon)|$, $|B_N(-1/4,0)|$ and $\mbf S_A^{(12)}$ in (\ref{eq:SAfinal}) with the choice of $\alpha_n$s as in (\ref{eq:alpha}).
\subsection{The estimation of $|B_N(*,*)|$}
\subsubsection{Estimation of $|B_N(0,0)|$:}
Since the $\alpha_n$s are supported on primes and prime squares we can write $|B_N(0,0)|$ as
\begin{equation*}
|B_N(0,0)|\le\sum_{p_1,p_2\le N^{1/2}}\frac{1}{(p_1p_2)^{3/2}}\mbf B_{p_1,p_2}+2\sum_{p_1,p_2\le N^{1/2}}\frac{1}{(p_1)^{3/2}p_2^3}\mbf B_{p_1,p_2^2}+\sum_{p_1,p_2\le N^{1/2}}\frac{1}{(p_1p_2)^3}\mbf B_{p_1^2,p_2^2},
\end{equation*}
where
\begin{equation*}
\mbf B_{n_1,n_2}:=\sum_{d|(n_1,n_2)}d^3\underset{a_2b_2^2d^2|n_1n_2}{\sum_{a_2,b_2}}\mu(\frac{n_1n_2}{a_2b_2^2})\sigma_{0}(a_2^2)b_2^{2}.
\end{equation*}
We have $\mbf B_{p_1,p_2}\ll 1$ if $p_1\neq p_2$ and is $\ll p^2$, if $p_1=p_2=p$; $\mbf B_{p_1,p_2^2}\ll p_2^2$ and $\mbf B_{p_1^2,p_2^2}\ll p_1^2p_2^2$. Thus
\begin{align*}
|B_N(0,0)|&\ll\sum_{p_1,p_2\le N^{1/2}}\frac{1}{(p_1p_2)^{3/2}}+\sum_{p\le N^{1/2}}\frac{1}{p}+2\sum_{p_1,p_2\le N^{1/2}}\frac{1}{(p_1)^{3/2}p_2}+\sum_{p_1,p_2\le N^{1/2}}\frac{1}{(p_1p_2)}\\
&\ll \log\log N.
\end{align*}
Here we use that $\sum_{p\le x}p^{-1}\asymp \log\log x$.
\subsubsection{Estimation of $\mid B_N(-\delta, 0) \mid$:} We have
\begin{equation*}
\mid B_N(-\delta, 0) \mid\le \sum_{n_1,n_2\le N}\frac{|\alpha_{n_1}\overline{\alpha_{n_2}}|}{(n_1n_2)^{3/2-2\delta}} \mbf B_{n_1,n_2},
\end{equation*}
where $\mbf B_{n_1,n_2}$ is as in the estimation for $|B_N(0,0)|$. Now evaluating similarly as in the case of $|B_N(0,0)|$, we have
\begin{align*}
\mid B_N(-\delta, 0) \mid&\ll\sum_{p_1,p_2\le N^{1/2}}\tfrac{1}{(p_1p_2)^{3/2-2\delta}}+\sum_{p\le N^{1/2}}\tfrac{1}{p^{1-4\delta}}+2\sum_{p_1,p_2\le N^{1/2}}\tfrac{1}{(p_1)^{3/2-2\delta}p_2^{1-4\delta}}+\sum_{p_1,p_2\le N^{1/2}}\tfrac{1}{(p_1p_2)^{1-4\delta}}\\
&\ll N^{-\frac{1}{2}+4\delta}.
\end{align*}
Here we use the fact that, for $s\neq 1$, $\sum_{p\le x}\frac{1}{p^s}\le \sum_{n\le x}\frac{1}{n^{s}}\asymp x^{1-s}$.
\subsubsection{Estimation of $|B_N(-\epsilon/2,\epsilon)|$:} We have $|B_N(-\epsilon/2,\epsilon)|$ is
\begin{equation*}
\le\sum_{p_1,p_2\le N^{1/2}}\tfrac{1}{(p_1p_2)^{3/2-\epsilon}}\mbf B^\epsilon_{p_1,p_2}+2\sum_{p_1,p_2\le N^{1/2}}\tfrac{1}{(p_1)^{3/2-\epsilon}p_2^{3-2\epsilon}}\mbf B^\epsilon_{p_1,p_2^2}+\sum_{p_1,p_2\le N^{1/2}}\tfrac{1}{(p_1p_2)^{3-2\epsilon}}\mbf B^\epsilon_{p_1^2,p_2^2},
\end{equation*}
We have $\mbf B^\epsilon_{p_1,p_2}\le (p_1p_2)^{1-\epsilon}$, $\mbf B^\epsilon_{p_1,p_2^2}\le (p_1p_2^2)^{1-\epsilon}$ and $\mbf B^\epsilon_{p_1^2,p_2^2}\le (p_1^2p_2^2)^{1-\epsilon}$. Thus we have
\begin{equation*}
|B_N(-\epsilon/2,\epsilon)|\ll \log\log N.
\end{equation*}
\subsubsection{Estimation of $|B_N(-\frac{1}{4},0)|$:}We have
\begin{equation*}
|B_N(-\frac{1}{4},0)|\le\sum_{n_1,n_2\le N}\frac{|\alpha_{n_1}\overline{\alpha_{n_2}}|}{n_1n_2}\sum_{d|(n_1,n_2)}d^{2}\underset{a_2b_2^2d^2|n_1n_2}{\sum_{a_2,b_2}}|\mu(\frac{n_1n_2}{a_2b_2^2})|\sigma_{0}(a_2^2)b_2.
\end{equation*}
Following the similar calculations as in the case of $|B_N(0,0)|$, we find that
\begin{equation}
|B_N(-\frac{1}{4},0)|\ll \log\log N.\nonumber
\end{equation}
\subsection{Estimation of $\mbf S_A^{(12)}$: }From (\ref{eq:S12}) we have
\begin{equation*}
\mbf S_A^{(12)}\ll \left(\sum_{n_1,n_2\le N}\frac{|\alpha_{n_1}\overline{\alpha_{n_2}}|}{(n_1n_2)^{5/2}+\epsilon}\sum_{d|(n_1,n_2)}d^{5+\epsilon}\underset{a_2b_2^2d^2|n_1n_2}{\sum_{a_2,b_2}}\mid\mu(\frac{n_1n_2}{a_2b_2^2})\mid a_2^{2+\epsilon}b_2^{4+2\epsilon}\right)k^{-1/3+\delta_1+\epsilon}.
\end{equation*}
We have that the inside summation is 
\begin{equation*}
\ll \sum_{p_1,p_2\le N^{1/2}}\tfrac{1}{(p_1p_2)^{5/2+\epsilon}}\mbf B'_{p_1,p_2}+2\sum_{p_1,p_2\le N^{1/2}}\tfrac{1}{(p_1)^{5/2+\epsilon}p_2^{5+2\epsilon}}\mbf B'_{p_1,p_2^2}+\sum_{p_1,p_2\le N^{1/2}}\tfrac{1}{(p_1p_2)^{5+2\epsilon}}\mbf B'_{p_1^2,p_2^2}
\end{equation*}
and $\mbf B'_{p_1,p_2}\ll (p_1p_2)^{2+\epsilon}$, $\mbf B'_{p_1,p_2^2}\ll (p_1p_2^2)^{2+\epsilon}$ and $\mbf B'_{p_1^2,p_2^2}\ll (p_1^2p_2^2)^{2+\epsilon}$. Thus
\begin{equation*}\label{eq:s12final}
\mbf S_A^{(12)}\ll \log\log N\cdot k^{-1/3+\delta_1+\epsilon}.
\end{equation*}
\subsection{Completion of \thmref{th:main}}
Now substituting in (\ref{eq:amplifiedterm}) with $N=k^\eta$ and using the fact that $\big|\sum_{p\le N^{1/2}}1\big|^2\asymp N/(\log N)^2$ and $\sum_{n_1,n_2\le N}|\alpha_{n_1}\overline{\alpha_{n_2}}|\ll \frac{N}{(\log N)^2}$, we have for any $1/2<\delta<1$
\begin{equation}\label{eq:sg}
\frac{12}{2k-1}S_{g_0}\ll k^{-\eta+\epsilon}+k^{-\delta_3+\epsilon}+k^{-1+\delta_2+\epsilon}+k^{-\frac{3\eta}{2}+4\eta\delta-2\delta+\epsilon}+k^{\frac{20\eta}{3}-\frac{1}{3}+\delta_1+\epsilon}.
\end{equation}
We make the following choice: $\delta_3=\frac{\delta_1\delta_2}{62}-\epsilon$ and note that since $1/2<\delta<1$, the fourth term in (\ref{eq:sg}) is irrelevant. Then we equate all the exponents of $k$ in (\ref{eq:sg}). A simple calculation shows that $\delta_1=27/91$ gives the answer. Also for this choice of $\delta_1$ we get $\eta\approx \frac{1}{210}$, $\delta_2\approx\frac{209}{210}$ and $\frac{1}{209}>\delta_3>\frac{1}{210}$. Thus we have
\begin{equation}\label{eq:powersaving}
S_{g_0}\ll k^{1-\frac{1}{210}+\epsilon}.
\end{equation}\qed

\end{document}